\documentclass[11pt,a4paper,fullpage]{article}
\usepackage[english]{babel}
\usepackage[utf8]{inputenc}
\usepackage{fancyhdr}
\usepackage[numbers]{natbib}

\usepackage{graphicx}
\usepackage{amsmath}
\usepackage{extsizes}
\usepackage{relsize}
\usepackage{multicol}
\usepackage{ragged2e}
\usepackage{fontenc}
\usepackage{mathcomp}
\usepackage{latexsym}
\usepackage{amssymb}
\usepackage{times}
\newtheorem{Teorema}{Theorem}[section]

\newtheorem{Posledica}[Teorema]{Corollary}

\newtheorem{Lema}[Teorema]{Lemma}
\newtheorem{Primedba}[Teorema]{Remark}
\usepackage{dsfont}
\numberwithin{equation}{section}
\hyphenchar\font=-1
\begin{document}
	\title {Upper triangular operator matrices and stability of their various spectra}	

	\author{Nikola Sarajlija\footnote{corresponding author: Nikola Sarajlija; University of Novi Sad, Faculty of Sciences, Novi Sad 21000, Serbia; {\it e-mail}: {\tt nikola.sarajlija@dmi.uns.ac.rs}} \footnote{the author gratefully acknowledge the financial support of the Ministry of Science, Technological Development and Innovation of the Republic of Serbia (Grants No. 451-03-66/2024-03/200125 and 451-03-65/2024-03/200125)} }

\maketitle

\begin{abstract}
Denote by $T_n^d(A)$ an upper triangular operator matrix of dimension $n\in\mathds{N}$ whose diagonal entries $D_i,\ 1\leq i\leq n$, are known, and $A=(A_{ij})_{1\leq i<j\leq n}$ is an unknown tuple of operators. This article is aimed at investigation of defect spectrum $\mathcal{D}^{\sigma_*}=\bigcup\limits_{i=1}^n\sigma_*(D_i)\setminus\sigma_*(T_n^d(A))$ , where $\sigma_*$ is a spectrum corresponding to various types of invertibility: (left, right) invertibility, (left, right) Fredholm invertibility, left/right Weyl invertibility. We give characterizations for each of the previous types, and provide some sufficent conditions for the stability of certain spectrum (the case $\mathcal{D}^{\sigma_*}=\emptyset$). The results are proved for all matrix dimensions $n\geq2$, and they hold in arbitrary Hilbert spaces without assuming separability, thus generalizing results from \cite{WU}, \cite{WU2}. We also retrieve a result from \cite{BAI} in the case $n=2$, and we provide a precise form of the well known 'filling in holes' result from \cite{HAN}.
\end{abstract}
\textit{$2020$ Math. Subj. Class:} 47A08, 47A53, 47A55, 47A05, 47A10.


\section{Introduction. Basic concepts}

Throughout this text, let $\mathcal{H},\mathcal{K},\mathcal{H}_1,...,\mathcal{H}_n$ be arbitrary Hilbert spaces. Collection of linear and bounded operators from $\mathcal{H}$ to $\mathcal{K}$ is denoted by $\mathcal{B}(\mathcal{H},\mathcal{K})$, where $\mathcal{B}(\mathcal{H}):=\mathcal{B}(\mathcal{H},\mathcal{H})$, and if $T\in\mathcal{B}(\mathcal{H},\mathcal{K})$ its null and range space are denoted by $\mathcal{N}(T)$ and $\mathcal{R}(T)$, respectively. We know that $\mathcal{N}(T)$ is a closed subspace of $\mathcal{H}$, and we denote its orthogonal dimension by $\alpha(T)$. Similarly, $\mathcal{R}(T)$ is a subspace of $\mathcal{K}$, and we denote its orthogonal codimension by $\beta(T)$. If $\alpha(T)$ ($\beta(T)$) is not finite we write $\alpha(T)=\infty$ ($\beta(T)=\infty$). We define $\mathrm{ind}(T):=\alpha(T)-\beta(T)$ as the index of $T$, if at least one of the quantities $\alpha(T),\beta(T)$ is finite. Notice that $\mathrm{ind}(T)\in\mathds{Z}\cup\lbrace+\infty,-\infty\rbrace.$

Let $T\in\mathcal{B}(\mathcal{H})$. Size $\alpha(T)$ ($\beta(T)$) introduced in the previous paragraph measures how close $T$ is to a left (right) invertible operator. If $\alpha(T)=0$ and $\mathcal{R}(T)$ is closed, then $T$ is left invertible, write $T\in\mathcal{G}_l(\mathcal{H})$. If $\beta(T)=0$, then $T$ is right invertible, write $T\in\mathcal{G}_r(\mathcal{H})$. If $\alpha(T)=\beta(T)=0$ then $T$ is invertible, write $T\in\mathcal{G}(\mathcal{H})=\mathcal{G}_l(\mathcal{H})\cap\mathcal{G}_r(\mathcal{H})$. Sizes  $\alpha$ and $\beta$ are systematically studied in Fredholm theory, the basic notation of which we introduce right now \cite{ZANA}. Families of left and right Fredholm operators, respectively, are defined as
$$
\begin{aligned}
\Phi_l(\mathcal{H})=\lbrace T\in\mathcal{B}(\mathcal{H}):\ \alpha(T)<\infty\ and\ \mathcal{R}(T)\ is \ closed\rbrace
\end{aligned}
$$
and
$$
\begin{aligned}
\Phi_r(\mathcal{H})=\lbrace T\in\mathcal{B}(\mathcal{H}): \beta(T)<\infty\rbrace.
\end{aligned}
$$
The set of Fredholm operators is
$$\Phi(\mathcal{H})=\Phi_l(\mathcal{H})\cap\Phi_r(\mathcal{H})=\lbrace T\in\mathcal{B}(\mathcal{H}): \alpha(T)<\infty\ and\ \beta(T)<\infty\rbrace.$$
Families of left and right Weyl operators, respectively, are defined as
$$\Phi_l^-(\mathcal{H})=\lbrace T\in\Phi_l(\mathcal{H}):\ \mathrm{ind}(T)\leq0\rbrace$$
and
$$\Phi_r^+(\mathcal{H})=\lbrace T\in\Phi_r(\mathcal{H}):\ \mathrm{ind}(T)\geq0\rbrace.$$

Corresponding spectra of an operator $T\in\mathcal{B}(\mathcal{H})$ are defined as follows:\\
-the left spectrum: $\sigma_{l}(T)=\lbrace\lambda\in\mathds{C}: \lambda-T\not\in\mathcal{G}_{l}(\mathcal{H})\rbrace$;\\
-the right spectrum: $\sigma_{r}(T)=\lbrace\lambda\in\mathds{C}: \lambda-T\not\in\mathcal{G}_{r}(\mathcal{H})\rbrace$;\\
-the spectrum: $\sigma(T)=\lbrace\lambda\in\mathds{C}: \lambda-T\not\in\mathcal{G}(\mathcal{H})\rbrace$;\\
-the left essential spectrum: $\sigma_{le}(T)=\lbrace\lambda\in\mathds{C}: \lambda-T\not\in\Phi_{l}(\mathcal{H})\rbrace$;\\
-the right essential spectrum: $\sigma_{re}(T)=\lbrace\lambda\in\mathds{C}: \lambda-T\not\in\Phi_{r}(\mathcal{H})\rbrace$;\\
-the essential spectrum: $\sigma_{e}(T)=\lbrace\lambda\in\mathds{C}: \lambda-T\not\in\Phi(\mathcal{H})\rbrace$;\\
-the left Weyl spectrum: $\sigma_{lw}(T)=\lbrace\lambda\in\mathds{C}: \lambda-T\not\in\Phi_l^-(\mathcal{H})\rbrace$;\\
-the right Weyl spectrum: $\sigma_{rw}(T)=\lbrace\lambda\in\mathds{C}: \lambda-T\not\in\Phi_r^+(\mathcal{H})\rbrace$;\\
All of these spectra are compact nonempty subsets of the complex plane.\\
Their complements are denoted by $\rho_l(T),\rho_r(T), \rho(T), \rho_{le}(T), \rho_{re}(T), \rho_e(T), \rho_{lw}(T), \rho_{rw}(T)$, respectively.

In this article we investigate spectral properties of upper triangular block operators whose diagonal entries are known while the others are not. Denote by $T_n^d(A)$ such an operator, where subscript $n$ denotes the matrix dimension and $A=(A_{ij})_{1\leq i<j\leq n}$ is a tuple of unknown operators above the main diagonal. In other words,
\begin{equation}\label{OSNOVNA}
T_n^d(A)=
\begin{bmatrix} 
    D_1 & A_{12} & A_{13} & ... & A_{1,n-1} & A_{1n}\\
    0 & D_2 & A_{23} & ... & A_{2,n-1} & A_{2n}\\
    0 &  0 & D_3 & ... & A_{3,n-1} & A_{3n}\\
    \vdots & \vdots & \vdots & \ddots & \vdots & \vdots\\
    0 & 0 & 0 & ... & D_{n-1} & A_{n-1,n}\\
    0 & 0 & 0 & ... & 0 & D_n      
\end{bmatrix}\in\mathcal{B}(\mathcal{H}_1\oplus \mathcal{H}_2\oplus\cdots\oplus \mathcal{H}_n),
\end{equation}
where it is understood that $A_{ij}\in\mathcal{B}(\mathcal{H}_j,\mathcal{H}_i)$, $1\leq i<j\leq n$.
For the convenience, we denote by $\mathcal{B}_n$ the collection of all described tuples $A=(A_{ij})$. Notice that this is the same notation already used by the present author in \cite{SARAJLIJA2}, \cite{SARAJLIJA3}, \cite{SARAJLIJA4}.

Investigation of spectral properties of block operators began with $2\times2$ block operators (see \cite{OPERATORTHEORY}, \cite{DU}, \cite{HAN} and many others).  General $n\times n$ block operators were first studied by Benhida et al. in 2005 \cite{BENHIDA}, but  afterwards they were mainly neglected until a few years ago (\cite{HUANG}, \cite{WU}, \cite{WU2}). This article has aim to complete investigation of basic spectral properties of $T_n^d(A)$ that appear in Fredholm theory.  Such an investigation has been previously conducted by the author of this article in \cite{SARAJLIJA2}, \cite{SARAJLIJA3}. Our main task will be to study defect spectrum $\mathcal{D}^{\sigma_*}=\bigcup\limits_{i=1}^n\sigma_*(D_i)\setminus\sigma_*(T_n^d(A))$, where $\sigma_*$ is one of the spectra introduced above. In the case $\mathcal{D}^{\sigma_*}=\emptyset$ we say that spectrum $\sigma_*$ is stable \cite{ZGUITTI2}. Stability of spectrum has usually been explored in a connection with the single valued extension property or SVEP for short (\cite{SLAVISA}, \cite{SLAVISA2}, \cite{ZGUITTI2}). Here we use a different method. It turns out that a relation between $\bigcup_{k=1}^n\sigma_*(D_k)$ and $\sigma_*(T_n^d(A))$ can be deduced using some characterization statements expressed for appropriate types of invertibility. Thus, we exploit results from \cite{SARAJLIJA2}, \cite{SARAJLIJA3} of the present author in order to gain information about defect spectrum.

Article is organized as follows. We end this section with a useful lemma that is used in the sequel. Its proof for $n=2$ that is conducted in \cite{OPERATORTHEORY} extends to an arbitrary dimension $n\geq2$ as well. Afterwards, in the next section we give results related to left/right Weyl, (left,right) Fredholm and (left,right) defect spectrum, respectively.
\begin{Lema}\label{POMOCNALEMA}
Let $T_n^d(A)\in\mathcal{B}(\mathcal{H}_1\oplus\cdots\oplus \mathcal{H}_n).$ Then:
\begin{itemize}
\item[(i)] $\sigma_{le}(D_1)\subseteq\sigma_{le}(T_n^d(A))\subseteq\bigcup\limits_{k=1}^n\sigma_{le}(D_k)$;
\item[(ii)] $\sigma_{re}(D_n)\subseteq\sigma_{re}(T_n^d(A))\subseteq\bigcup\limits_{k=1}^n\sigma_{re}(D_k);$
\item[(iii)]$\sigma_{le}(D_1)\cup\sigma_{re}(D_n)\subseteq\sigma_{e}(T_n^d(A))\subseteq\bigcup\limits_{k=1}^n\sigma_{e}(D_k)$;
\item[(iv)]$\sigma_{lw}(T_n^d(A))\subseteq\bigcup\limits_{k=1}^n\sigma_{lw}(D_k)$;
\item[(v)]$\sigma_{rw}(T_n^d(A))\subseteq\bigcup\limits_{k=1}^n\sigma_{rw}(D_k)$;
\item[(vi)]$\sigma_{l}(D_1)\subseteq\sigma_{l}(T_n^d(A))\subseteq\bigcup\limits_{k=1}^n\sigma_{l}(D_k)$;
\item[(vii)]$\sigma_{r}(D_n)\subseteq\sigma_{r}(T_n^d(A))\subseteq\bigcup\limits_{k=1}^n\sigma_{r}(D_k)$;
\item[(viii)]$\sigma_l(D_1)\cup\sigma_r(D_n)\subseteq\sigma(T_n^d(A))\subseteq\bigcup\limits_{k=1}^n\sigma(D_k)$.
\end{itemize}
\end{Lema}
 
\section{Defect spectra and stability results}

In this section, $\mathcal{H}_1,...,\mathcal{H}_n$ are infinite dimensional Hilbert spaces. Occasionally, we will need an assumption that the former are separable, in which case we shall emphasize this fact.

\subsection{Left and right Weyl spectrum}
In this subsection we generalize results from \cite[Section 3]{WU2} to arbitrary Hilbert spaces. We report that Corollaries 3.3 and 3.8 in \cite{WU2} do not hold with the equivalence: 'only if' part is not valid. The reason for this is that the proofs of these corollaries summon \cite[Theorems 2.5, 2.6]{WU2} which do not hold with an equality (see \cite[Corollaries 2.3, 2.10]{SARAJLIJA2} for corrected versions). Corollaries 4, 8 and 12 from \cite[Section 3]{WU} are not valid for analogous reasons. In the sequel we provide correct forms of these statements.

\begin{Teorema}(\cite[Corollary 2.3]{SARAJLIJA2})\label{POSLEDICA'}
Let $D_1\in\mathcal{B}(\mathcal{H}_1),\ D_2\in\mathcal{B}(\mathcal{H}_2),...,D_n\in\mathcal{B}(\mathcal{H}_n)$. Then
\begin{equation}\label{STARIJE1}
\sigma_{le}(D_1)\cup\Big(\bigcup\limits_{k=2}^{n+1}\delta_k\Big)\subseteq\\\bigcap\limits_{A\in\mathcal{B}_n}\sigma_{lw}(T_n^d(A)),
\end{equation}
where
$$
\delta_k:=\Big\lbrace\lambda\in\mathds{C}:\ \alpha(D_k-\lambda)=\infty\ and\  \sum\limits_{s=1}^{k-1}\beta(D_s-\lambda)<\infty\Big\rbrace,\ 2\leq k\leq n,
$$
$$
\delta_{n+1}:=\Big\lbrace\lambda\in\mathds{C}:\ \sum\limits_{s=1}^n\beta(D_s-\lambda)<\sum\limits_{s=1}^{n}\alpha(D_s-\lambda)\Big\rbrace.
$$
\end{Teorema}

\begin{Teorema}\label{pOSLEDICA'}
Let $D_1\in\mathcal{B}(\mathcal{H}_1),\ D_2\in\mathcal{B}(\mathcal{H}_2),...,D_n\in\mathcal{B}(\mathcal{H}_n)$. If $\mathcal{H}_1,...,\mathcal{H}_n$ are separable and $\mathcal{R}(D_s-\lambda)$ are closed for $2\leq s\leq n$, $\lambda\in\mathds{C},$ then
\begin{equation}\label{sTARIJE1}
\sigma_{le}(D_1)\cup\Big(\bigcup\limits_{k=2}^{n+1}\delta_k\Big)=\\\bigcap\limits_{A\in\mathcal{B}_n}\sigma_{lw}(T_n^d(A)),
\end{equation}
where $\delta_k$, $2\leq k\leq n+1$, are defined as in Theorem \ref{POSLEDICA'}.
\end{Teorema}
\textbf{Proof. }This is obvious from \cite[Corollary 2.5]{SARAJLIJA3}. $\square$

Now we are able to prove the following generalization to arbitrary Hilbert spaces of \cite[Theorem 3.1]{WU2}.

\begin{Teorema}\label{PRVA}
Let $D_1\in\mathcal{B}(\mathcal{H}_1),\ D_2\in\mathcal{B}(\mathcal{H}_2),...,D_n\in\mathcal{B}(\mathcal{H}_n)$. Then
\begin{equation}\label{PRVAJEDAN}
\bigcup\limits_{k=1}^n\sigma_{lw}(D_k)=\sigma_{lw}(T_n^d(A))\cup\Delta_1\cup\Delta_2
\end{equation}
holds for every $A\in\mathcal{B}_n$, where 
$$
\begin{aligned}
\Delta_1=\bigcup\limits_{k=2}^{n}\{\lambda\in\mathds{C}:\ \alpha(D_k-\lambda)=\infty,\ \alpha(D_s-\lambda)<\infty\ for\ 2\leq s\leq k-1\ and\\ \sum\limits_{s=1}^{k-1}\beta(D_s-\lambda)=\infty\}\cap\rho_{le}(D_1)\cap\{\lambda\in\mathds{C}:\ \sum\limits_{s=1}^n\beta(D_s-\lambda)\geq\sum\limits_{s=1}^{n}\alpha(D_s-\lambda)\},
\end{aligned}
$$
$$
\begin{aligned}
\Delta_2=\bigcup\limits_{k=1}^{n}\Big\{\lambda\in\mathds{C}:\ \alpha(D_s-\lambda)<\infty\ for\ all\ 1\leq s\leq n, \sum\limits_{s=1}^n\beta(D_s-\lambda)\geq\\\sum\limits_{s=1}^{n}\alpha(D_s-\lambda)\ and\ \Big(\alpha(D_k-\lambda)>\beta(D_k-\lambda)\ or\ \mathcal{R}(D_k-\lambda)\ is\ not\ closed\Big)\Big\}.
\end{aligned}
$$
\end{Teorema}
\begin{Primedba}
Condition '$\alpha(D_s-\lambda)<\infty\ for\ 2\leq s\leq k-1$' in $\Delta_1$ is omitted when $k=2$.
\end{Primedba}

\textbf{Proof.} First of all, notice that $\sigma_{lw}(T_n^d(A))\subseteq\bigcup\limits_{k=1}^n\sigma_{lw}(D_k)$ according to Lemma \ref{POMOCNALEMA}. Next, we prove that $\Delta_1\cup\Delta_2\subseteq\bigcup\limits_{k=1}^n\sigma_{lw}(D_k)$. This is, however, obvious from the definition of a left Weyl operator. Namely, if $\lambda\in\Delta_1$ then $\alpha(D_k-\lambda)=\infty$ for some $k\in\lbrace1,...,n\rbrace,$ which means that $D_k-\lambda$ is not left Weyl, thus $\lambda\in\sigma_{lw}(D_k)\subseteq\bigcup\limits_{k=1}^n\sigma_{lw}(D_k)$. Similarly, $\lambda\in\Delta_2$ implies ($\alpha(D_k-\lambda)>\beta(D_k-\lambda)$ or $\mathcal{R}(D_k-\lambda)$ is not closed) for some $k\in\lbrace1,...,n\rbrace$ which means that $D_k-\lambda$ is not left Weyl, thus $\lambda\in\sigma_{lw}(D_k)\subseteq\bigcup\limits_{k=1}^n\sigma_{lw}(D_k)$ again. Now, let us prove the opposite inclusion in (\ref{PRVAJEDAN}). Assume that $\lambda\in\bigcup\limits_{k=1}^n\sigma_{lw}(D_k)\setminus\sigma_{lw}(T_n^d(A))$. We want to prove that in this case, $\lambda\in\Delta_1\cup\Delta_2$. However, by Theorem \ref{POSLEDICA'}, we get that $\lambda$ does not belong to the left side of (\ref{STARIJE1}), which, using the notation from (\ref{STARIJE1}), means that $\lambda\in\rho_{le}(D_1)\cap\bigcap\limits_{k=2}^{n+1}\delta_k^c$. This,  together with our assumption $\lambda\in\bigcup\limits_{k=1}^n\sigma_{lw}(D_k)$, easily gives $\lambda\in\Delta_1\cup\Delta_2$. $\square$

Now, we can give a sufficient condition for the stability of the left Weyl spectrum.

\begin{Posledica}
Let $D_1\in\mathcal{B}(\mathcal{H}_1),\ D_2\in\mathcal{B}(\mathcal{H}_2),...,D_n\in\mathcal{B}(\mathcal{H}_n)$. Then
$$
\bigcup\limits_{k=1}^n\sigma_{lw}(D_k)=\sigma_{lw}(T_n^d(A))
$$
holds for every $A\in\mathcal{B}_n$ if
$$
\Delta_1\cup\Delta_2=\emptyset,
$$
where $\Delta_1,\Delta_2$ are defined as in Theorem \ref{PRVA}.
\end{Posledica}

If we summon the separability condition, then we are able to state the following.

\begin{Posledica}(\cite[Theorem 3.1]{WU2}, corrected version)
Let $D_1\in\mathcal{B}(\mathcal{H}_1),\ D_2\in\mathcal{B}(\mathcal{H}_2),...,D_n\in\mathcal{B}(\mathcal{H}_n)$. If $\mathcal{H}_1,...,\mathcal{H}_n$ are separable and $\mathcal{R}(D_s-\lambda),\ 2\leq s\leq n,\ \lambda\in\mathds{C}$ are closed, then
$$
\bigcup\limits_{k=1}^n\sigma_{lw}(D_k)=\sigma_{lw}(T_n^d(A))
$$
holds for every $A\in\mathcal{B}_n$ if and only if
$$
\Delta_1\cup\Delta_2=\emptyset,
$$
where $\Delta_1,\Delta_2$ are defined as in Theorem \ref{PRVA}.
\end{Posledica}
\textbf{Proof. }Sufficiency is clear, and necessity follows from Theorem \ref{PRVA} and Theorem \ref{pOSLEDICA'}. $\square$

By duality, we obtain results related to the stability of the right Weyl spectrum. We begin with the following generalization of \cite[Theorem 3.6]{WU2} to arbitrary Hilbert spaces.

\begin{Teorema}\label{DRUGA}
Let $D_1\in\mathcal{B}(\mathcal{H}_1),\ D_2\in\mathcal{B}(\mathcal{H}_2),...,D_n\in\mathcal{B}(\mathcal{H}_n)$. Then
$$
\bigcup\limits_{k=1}^n\sigma_{rw}(D_k)=\sigma_{rw}(T_n^d(A))\cup\Delta_1\cup\Delta_2
$$
holds for every $A\in\mathcal{B}_n$, where 
$$
\begin{aligned}
\Delta_1=\bigcup\limits_{k=1}^{n-1}\{\lambda\in\mathds{C}:\ \beta(D_k-\lambda)=\infty,\ \beta(D_s-\lambda)<\infty\ for\ k+1\leq s\leq n-1\ and\\ \sum\limits_{s=k+1}^{n}\alpha(D_s-\lambda)=\infty\}\cap\rho_{re}(D_n)\cap\{\lambda\in\mathds{C}:\ \sum\limits_{s=1}^n\alpha(D_s-\lambda)\geq\sum\limits_{s=1}^{n}\beta(D_s-\lambda)\},
\end{aligned}
$$
$$
\begin{aligned}
\Delta_2=\bigcup\limits_{k=1}^{n}\{\lambda\in\mathds{C}:\ \beta(D_s-\lambda)<\infty\ for\ all\ 1\leq s\leq n,\ \sum\limits_{s=1}^n\alpha(D_s-\lambda)\geq\\\sum\limits_{s=1}^{n}\beta(D_s-\lambda)\ and\ \beta(D_k-\lambda)>\alpha(D_k-\lambda)\}.
\end{aligned}
$$
\end{Teorema}
\begin{Primedba}
Condition '$\beta(D_s-\lambda)<\infty\ for\ k+1\leq s\leq n-1$' in $\Delta_1$ is omitted when $k=n-1$.
\end{Primedba}
Now, we can give a sufficient condition for the stability of the right Weyl spectrum.

\begin{Posledica}
Let $D_1\in\mathcal{B}(\mathcal{H}_1),\ D_2\in\mathcal{B}(\mathcal{H}_2),...,D_n\in\mathcal{B}(\mathcal{H}_n)$. Then
$$
\bigcup\limits_{k=1}^n\sigma_{rw}(D_k)=\sigma_{rw}(T_n^d(A))
$$
holds for every $A\in\mathcal{B}_n$ if 
$$
\Delta_1\cup\Delta_2=\emptyset,
$$
where $\Delta_1,\Delta_2$ are defined as in Theorem \ref{DRUGA}.
\end{Posledica}

If we include the separability assumption, we obtain characterization for the stability of the right Weyl spectrum.

\begin{Posledica}(\cite[Theorem 3.6]{WU2}, corrected version)
Let $D_1\in\mathcal{B}(\mathcal{H}_1),\ D_2\in\mathcal{B}(\mathcal{H}_2),...,D_n\in\mathcal{B}(\mathcal{H}_n)$. If $\mathcal{H}_1,...,\mathcal{H}_n$ are separable and $\mathcal{R}(D_s-\lambda),\ 1\leq s\leq n-1,\ \lambda\in\mathds{C}$ are closed, then
$$
\bigcup\limits_{k=1}^n\sigma_{rw}(D_k)=\sigma_{rw}(T_n^d(A))
$$
holds for every $A\in\mathcal{B}_n$ if and only if
$$
\Delta_1\cup\Delta_2=\emptyset,
$$
where $\Delta_1,\Delta_2$ are defined as in Theorem \ref{DRUGA}.
\end{Posledica}

\subsection{Fredholm spectrum}

In this subsection we generalize results from \cite[Section 3]{WU} to arbitrary Hilbert spaces.  We prove statements related to left Fredholm invertibility, and then by duality obtain corresponding statements related to right Fredholm invertibility. Finally, we finish this subsection with investigation of the essential spectra.

We start with the following two known results.

\begin{Teorema}(\cite[Corollary 3.3]{SARAJLIJA2})\label{POSLEDICA5''}
Let $D_1\in\mathcal{B}(\mathcal{H}_1),\ D_2\in\mathcal{B}(\mathcal{H}_2),...,D_n\in\mathcal{B}(\mathcal{H}_n)$. Then
\begin{equation}\label{STARIJE2}
\sigma_{le}(D_1)\cup\Big(\bigcup\limits_{k=2}^{n}\delta_k\Big)\subseteq\bigcap\limits_{A\in\mathcal{B}_n}\sigma_{le}(T_n^d(A)),
\end{equation}
where
$$
\delta_k:=\Big\lbrace\lambda\in\mathds{C}:\ \alpha(D_k-\lambda)=\infty\ and\  \sum\limits_{s=1}^{k-1}\beta(D_s-\lambda)<\infty\Big\rbrace,\ 2\leq k\leq n.
$$
\end{Teorema}

\begin{Teorema}\label{pOSLEDICA5''}
Let $D_1\in\mathcal{B}(\mathcal{H}_1),\ D_2\in\mathcal{B}(\mathcal{H}_2),...,D_n\in\mathcal{B}(\mathcal{H}_n)$. If $\mathcal{H}_1,...,\mathcal{H}_n$ are separable and $\mathcal{R}(D_s-\lambda)$, $2\leq s\leq n$, $\lambda\in\mathds{C}$ are closed, then
\begin{equation}\label{STARIJE2}
\sigma_{le}(D_1)\cup\Big(\bigcup\limits_{k=2}^{n}\delta_k\Big)=\bigcap\limits_{A\in\mathcal{B}_n}\sigma_{le}(T_n^d(A)),
\end{equation}
where $\delta_k$, $2\leq k\leq n$, are defined as in Theorem \ref{POSLEDICA5''}.
\end{Teorema}
\textbf{Proof.} This is obvious from \cite[Corollary 2.20]{SARAJLIJA3}. $\square$

Now, we generalize \cite[Theorem 4]{WU} to arbitrary Hilbert spaces.

\begin{Teorema}\label{TRECA}
Let $D_1\in\mathcal{B}(\mathcal{H}_1),\ D_2\in\mathcal{B}(\mathcal{H}_2),...,D_n\in\mathcal{B}(\mathcal{H}_n)$. Then
\begin{equation}\label{TRECAJEDAN}
\bigcup\limits_{k=1}^n\sigma_{le}(D_k)=\sigma_{le}(T_n^d(A))\cup\Delta_1\cup\Delta_2
\end{equation}
holds for every $A\in\mathcal{B}_n$, where 
$$
\begin{aligned}
\Delta_1=\bigcup\limits_{k=2}^{n}\{\lambda\in\mathds{C}:\ \alpha(D_k-\lambda)=\infty,\ \alpha(D_s-\lambda)<\infty\ for\ 2\leq s\leq k-1\\ and\ \sum\limits_{s=1}^{k-1}\beta(D_s-\lambda)=\infty\}\cap\rho_{le}(D_1),
\end{aligned}
$$
$$
\begin{aligned}
\Delta_2=\bigcup\limits_{k=1}^{n}\{\lambda\in\mathds{C}:\ \alpha(D_s-\lambda)<\infty\ for\ all\ s=1,...,n\ and\\ \mathcal{R}(D_k-\lambda)\ is\ not\ closed\}.
\end{aligned}
$$
\end{Teorema}

\begin{Primedba}
Condition '$\alpha(D_s-\lambda)<\infty\ for\ 2\leq s\leq k-1$' in $\Delta_1$ is omitted when $k=2$.
\end{Primedba}

\textbf{Proof.} Obviously, $\Delta_1\cup\Delta_2\subseteq\bigcup\limits_{k=1}^n\sigma_{le}(D_k)$, and $\sigma_{le}(T_n^d(A))\subseteq\bigcup\limits_{k=1}^n\sigma_{le}(D_k)$ according to Lemma \ref{POMOCNALEMA}. Assume that $\lambda\in\bigcup\limits_{k=1}^n\sigma_{le}(D_k)\setminus\sigma_{le}(T_n^d(A))$. Then by Theorem \ref{POSLEDICA5''} we get that $\lambda$ does not belong to the left side of (\ref{STARIJE2}), which together with observation $\lambda\in\bigcup\limits_{k=1}^n\sigma_{le}(D_k)$ easily gives $\lambda\in\Delta_1\cup\Delta_2$. $\square$

Now, we can give sufficient condition for the stability of the left Fredholm spectrum.

\begin{Posledica}
Let $D_1\in\mathcal{B}(\mathcal{H}_1),\ D_2\in\mathcal{B}(\mathcal{H}_2),...,D_n\in\mathcal{B}(\mathcal{H}_n)$. Then
$$
\bigcup\limits_{k=1}^n\sigma_{le}(D_k)=\sigma_{le}(T_n^d(A))
$$
holds for every $A\in\mathcal{B}_n$ if
$$
\Delta_1\cup\Delta_2=\emptyset,
$$
where $\Delta_1,\Delta_2$ are defined as in Theorem \ref{TRECA}.
\end{Posledica}

\begin{Posledica}(\cite[Corollary 4]{WU}, corrected version)
Let $D_1\in\mathcal{B}(\mathcal{H}_1),\ D_2\in\mathcal{B}(\mathcal{H}_2),...,D_n\in\mathcal{B}(\mathcal{H}_n)$. If $\mathcal{H}_1,...,\mathcal{H}_n$ are separable and $\mathcal{R}(D_s-\lambda),\ 2\leq s\leq n,\ \lambda\in\mathds{C}$ are closed, then
$$
\bigcup\limits_{k=1}^n\sigma_{le}(D_k)=\sigma_{le}(T_n^d(A))
$$
holds for every $A\in\mathcal{B}_n$ if and only if
$$
\Delta_1\cup\Delta_2=\emptyset,
$$
where $\Delta_1,\Delta_2$ are defined as in Theorem \ref{TRECA}.
\end{Posledica}
\textbf{Proof.} Sufficiency is obvious, and necessity follows from Theorem \ref{TRECA} and Theorem \ref{pOSLEDICA5''}. $\square$ 

We provide the following results for the right Fredholm spectrum. First we generalize \cite[Theorem 5]{WU} to arbitrary Hilbert spaces.

\begin{Teorema}\label{CETVRTA}
Let $D_1\in\mathcal{B}(\mathcal{H}_1),\ D_2\in\mathcal{B}(\mathcal{H}_2),...,D_n\in\mathcal{B}(\mathcal{H}_n)$. Then
$$
\bigcup\limits_{k=1}^n\sigma_{re}(D_k)=\sigma_{re}(T_n^d(A))\cup\Delta_1\cup\Delta_2
$$
holds for every $A\in\mathcal{B}_n$, where 
$$
\begin{aligned}
\Delta_1=\bigcup\limits_{k=1}^{n-1}\{\lambda\in\mathds{C}:\ \beta(D_k-\lambda)=\infty,\ \beta(D_s-\lambda)<\infty\ for\ k+1\leq s\leq n-1\ and\\ \sum\limits_{s=k+1}^{n}\alpha(D_s-\lambda)=\infty\}\cap\rho_{re}(D_n),
\end{aligned}
$$
$$
\begin{aligned}
\Delta_2=\bigcup\limits_{k=1}^n\{\lambda\in\mathds{C}:\ \beta(D_s-\lambda)<\infty\ for\ all\ 1\leq s\leq n\ and\ \mathcal{R}(D_k-\lambda)\ is\ not\ closed\}.
\end{aligned}
$$
\end{Teorema}
\begin{Primedba}
Condition '$\beta(D_s-\lambda)<\infty\ for\ k+1\leq s\leq n-1$' in $\Delta_1$ is omitted when $k=n-1$.
\end{Primedba}

Sufficient condition for the stability of the right Fredholm spectrum follows.

\begin{Posledica}
Let $D_1\in\mathcal{B}(\mathcal{H}_1),\ D_2\in\mathcal{B}(\mathcal{H}_2),...,D_n\in\mathcal{B}(\mathcal{H}_n)$. Then
$$
\bigcup\limits_{k=1}^n\sigma_{re}(D_k)=\sigma_{re}(T_n^d(A))
$$
holds for every $A\in\mathcal{B}_n$ if 
$$
\Delta_1\cup\Delta_2=\emptyset,
$$
where $\Delta_1,\Delta_2$ are defined as in Theorem \ref{CETVRTA}.
\end{Posledica}

Let us summon separability next.

\begin{Posledica}(\cite[Corollary 8]{WU}, corrected version)
Let $D_1\in\mathcal{B}(\mathcal{H}_1),\ D_2\in\mathcal{B}(\mathcal{H}_2),...,D_n\in\mathcal{B}(\mathcal{H}_n)$. Assume that $\mathcal{H}_1,...,\mathcal{H}_n$ are separable and $\mathcal{R}(D_s-\lambda),\ 1\leq s\leq n-1,\ \lambda\in\mathds{C}$ are closed. Then
$$
\bigcup\limits_{k=1}^n\sigma_{re}(D_k)=\sigma_{re}(T_n^d(A))
$$
holds for every $A\in\mathcal{B}_n$ if and only if
$$
\Delta_1\cup\Delta_2=\emptyset,
$$
where $\Delta_1,\Delta_2$ are defined as in Theorem \ref{CETVRTA}.
\end{Posledica}

To end this section, we provide statements dealing with the essential spectrum. We begin with

\begin{Teorema}(\cite[Corollary 3.17]{SARAJLIJA2})\label{POSLEDICA9''}
Let $D_1\in\mathcal{B}(\mathcal{H}_1),\ D_2\in\mathcal{B}(\mathcal{H}_2),...,D_n\in\mathcal{B}(\mathcal{H}_n)$. Then
\begin{equation}\label{STARIJE3}
\sigma_{le}(D_1)\cup\sigma_{re}(D_n)\cup\Big(\bigcup\limits_{k=2}^{n-1}\delta_k\Big)\cup\delta_n\subseteq\\\bigcap_{A\in\mathcal{B}_n}\sigma_e(T_n^d(A))
\end{equation}
where
$$
\begin{aligned}
\delta_k=\Big\lbrace\lambda\in\mathds{C}:\ \alpha(D_k-\lambda)=\infty\ and\ \sum_{s=1}^{k-1}\beta(D_s-\lambda)<\infty\Big\rbrace\cup\\
\Big\lbrace\lambda\in\mathds{C}:\ \beta(D_k-\lambda)=\infty\ and\ \sum_{s=k+1}^n\alpha(D_s-\lambda)<\infty\Big\rbrace,\quad 2\leq k\leq n-1,
\end{aligned}
$$
$$
\begin{aligned}
\delta_n=\Big\lbrace\lambda\in\mathds{C}:\ \alpha(D_n-\lambda)=\infty\ and\ \sum_{s=1}^{n-1}\beta(D_s-\lambda)<\infty\Big\rbrace\cup\\\Big\lbrace\lambda\in\mathds{C}:\ \beta(D_1-\lambda)=\infty\ and\ \sum_{s=2}^n\alpha(D_s-\lambda)<\infty\Big\rbrace,
\end{aligned}
$$
\end{Teorema}

\begin{Teorema}\label{pOSLEDICA9''}
Let $D_1\in\mathcal{B}(\mathcal{H}_1),\ D_2\in\mathcal{B}(\mathcal{H}_2),...,D_n\in\mathcal{B}(\mathcal{H}_n)$. If $\mathcal{H}_1,...,\mathcal{H}_n$ are separable and $\mathcal{R}(D_s-\lambda)$, $2\leq s\leq n-1$, $\lambda\in\mathds{C}$, are closed, then
\begin{equation}\label{STARIJE3}
\sigma_{le}(D_1)\cup\sigma_{re}(D_n)\cup\Big(\bigcup\limits_{k=2}^{n-1}\delta_k\Big)\cup\delta_n=\\ \bigcap_{A\in\mathcal{B}_n}\sigma_e(T_n^d(A)),
\end{equation}
where $\delta_k$, $2\leq k\leq n$, are defined as in Theorem \ref{POSLEDICA9''}.
\end{Teorema}
\textbf{Proof. }This is obvious from \cite[Corollary 2.34]{SARAJLIJA3}. $\square$

\begin{Teorema}\label{PETA}
Let $D_1\in\mathcal{B}(\mathcal{H}_1),\ D_2\in\mathcal{B}(\mathcal{H}_2),...,D_n\in\mathcal{B}(\mathcal{H}_n)$. Then
$$
\bigcup\limits_{k=1}^n\sigma_{e}(D_k)=\sigma_{e}(T_n^d(A))\cup\Delta
$$
holds for every $A\in\mathcal{B}_n$, where
$$
\Delta=(\Delta_1\cup\Delta_2)\cap\rho_{le}(D_1)\cap\rho_{re}(D_n),
$$ 
$$
\begin{aligned}
\Delta_1=\bigcup\limits_{k=2}^{n-1}\Big\{\lambda\in\mathds{C}:\ \Big(\alpha(D_k-\lambda)=\sum\limits_{s=1}^{k-1}\beta(D_s-\lambda)=\infty\ and\ \alpha(D_s-\lambda)<\infty\\ for\ 2\leq s\leq k-1\Big)\ or\ \Big(\beta(D_k-\lambda)=\sum\limits_{s=k+1}^{n}\alpha(D_s-\lambda)=\infty\ and\ \beta(D_s-\lambda)<\infty\\ for\ k+1\leq s\leq n-1\Big)\Big\},
\end{aligned}
$$
$$
\begin{aligned}
\Delta_2=\Big\{\lambda\in\mathds{C}:\ \Big(\alpha(D_n-\lambda)=\sum\limits_{s=1}^{n-1}\beta(D_s-\lambda)=\infty\ and\ \alpha(D_s-\lambda)<\infty\\ for\ 2\leq s\leq n-1\Big)\ or\ \Big(\beta(D_1-\lambda)=\sum\limits_{s=2}^{n}\alpha(D_s-\lambda)=\infty\  and\ \beta(D_s-\lambda)<\infty\\ for\ 2\leq s\leq n-1\Big)\Big\}.
\end{aligned}
$$
\end{Teorema}
\textbf{Proof.} First of all, $\sigma_e(T_n^d(A))\subseteq\bigcup\limits_{k=1}^n\sigma_{e}(D_k)$ according to Lemma \ref{POMOCNALEMA}. Next, we prove that $\Delta\subseteq\bigcup\limits_{k=1}^n\sigma_{e}(D_k)$. This is, however, obvious from the definition of a Fredholm operator. Namely, if $\lambda\in \Delta_1$ then $\alpha(D_k-\lambda)=\infty$ or $\beta(D_k-\lambda)=\infty$ for some $k\in\lbrace2,...,n-1\rbrace$ which means that $D_k-\lambda$ is not Fredholm, thus $\lambda\in\sigma_e(D_k)\subseteq\bigcup\limits_{k=1}^n\sigma_{e}(D_k)$. Similarly, $\lambda\in\Delta_2$ implies $\alpha(D_n-\lambda)=\infty$ or $\beta(D_1-\lambda)=\infty$ which means that either $\lambda\in\sigma_{le}(D_n)\subseteq\bigcup\limits_{k=1}^n\sigma_{e}(D_k)$ or $\lambda\in\sigma_{re}(D_1)\subseteq\bigcup\limits_{k=1}^n\sigma_{e}(D_k)$, respectively. Now, let us prove the opposite inclusion. Assume that $\lambda\in\bigcup\limits_{k=1}^n\sigma_{e}(D_k)\setminus\sigma_{e}(T_n^d(A))$. However, by Theorem \ref{POSLEDICA9''}, we get that $\lambda$ does not belong to the left side of (\ref{STARIJE3}), which, using the notation from (\ref{STARIJE3}), means that $\lambda\in\rho_{le}(D_1)\cap\rho_{re}(D_n)\cap\bigcap\limits_{k=2}^{n}\delta_k^c$. This, together with our assumption $\lambda\in\bigcup\limits_{k=1}^n\sigma_{e}(D_k)$, easily gives $\lambda\in\Delta$. $\square$

Now we can give sufficient condition for the stability of the essential spectrum.
\begin{Posledica}
Let $D_1\in\mathcal{B}(\mathcal{H}_1),\ D_2\in\mathcal{B}(\mathcal{H}_2),...,D_n\in\mathcal{B}(\mathcal{H}_n)$. Then
$$
\bigcup\limits_{k=1}^n\sigma_{e}(D_k)=\sigma_{e}(T_n^d(A))
$$
holds for every $A\in\mathcal{B}_n$ if 
$$
\Delta=\emptyset,
$$
where $\Delta$ is defined as in Theorem \ref{PETA}.
\end{Posledica}

\begin{Posledica}(\cite[Corollary 12]{WU}, corrected version)
Let $D_1\in\mathcal{B}(\mathcal{H}_1),\ D_2\in\mathcal{B}(\mathcal{H}_2),...,D_n\in\mathcal{B}(\mathcal{H}_n)$. Assume that $\mathcal{H}_1,...,\mathcal{H}_n$ are separable and $\mathcal{R}(D_s-\lambda),\ 2\leq s\leq n-1,\ \lambda\in\mathds{C}$ are closed. Then
$$
\bigcup\limits_{k=1}^n\sigma_{e}(D_k)=\sigma_{e}(T_n^d(A))
$$
holds for every $A\in\mathcal{B}_n$ if and only if
$$
\Delta=\emptyset,
$$
where $\Delta$ is defined as in Theorem \ref{PETA}.
\end{Posledica}
\textbf{Proof.} Sufficiency is obvious, and necessity follows from Theorem \ref{PETA} and Theorem \ref{pOSLEDICA9''}. $\square$ 

Statements related to the Fredholm spectrum become especially elegant when $n=2$.
\begin{Teorema}(\cite[Corollary 3.2]{BAI})
Let $D_1\in\mathcal{B}(\mathcal{H}_1),\ D_2\in\mathcal{B}(\mathcal{H}_2)$. Then
$$
\begin{aligned}
\sigma_{e}(D_1)\cup\sigma_e(D_2)=\sigma_{e}(T_2^d(A))\cup\Big(\lbrace\lambda\in\mathds{C}:\ \beta(D_1-\lambda)=\alpha(D_2-\lambda)=\infty\rbrace\cap\\ \rho_{le}(D_1)\cap\rho_{re}(D_2)\Big)
\end{aligned}
$$
holds for every $A\in\mathcal{B}_2$.
\end{Teorema}
\begin{Posledica}
Let $D_1\in\mathcal{B}(\mathcal{H}_1),\ D_2\in\mathcal{B}(\mathcal{H}_2)$. Then
$$
\sigma_{e}(D_1)\cup\sigma_e(D_2)=\sigma_{e}(T_2^d(A))
$$
holds for every $A\in\mathcal{B}_2$ if
$$
\lbrace\lambda\in\mathds{C}:\ \beta(D_1-\lambda)=\alpha(D_2-\lambda)=\infty\rbrace\cap\\ \rho_{le}(D_1)\cap\rho_{re}(D_2)=\emptyset.
$$
\end{Posledica}
Thus, we recover Remark 3 from \cite[Section 3]{WU}.
\subsection{Special classes of diagonal operators}

It is not hard to see that all corollaries of the present section summon the same assumption: diagonal operators $D_s$ must have closed range together with all their translates $D_s-\lambda$ for appropriate indices $s\in\mathds{N}$. One might wonder if there are any operators other than the trivial ones satisfying such an assumption. However, this assumption holds true for some relatively large classes of operators, as we shall prove in the sequel.

\begin{Teorema}
If $T\in\mathcal{B}(\mathcal{H})$ is compact, then $\mathcal{R}(T-\lambda)$ is closed for every $\lambda\in\mathds{C}\lbrace0\rbrace$.
\end{Teorema}
\textbf{Proof.} Assume that $T\in\mathcal{H}$ is a compact operator. Then $\frac{1}{\lambda}T$ is compact for every $\lambda\in\mathds{C}\setminus\lbrace0\rbrace$ as well. Therefore, by the Fredholm alternative, $\mathcal{R}(\frac{1}{\lambda}T-I)$ is closed, and thus $\mathcal{R}(T-\lambda)$ is closed as well. This finishes our proof. $\square$

Since every finite rank operator is compact, we have an immediate consequence.
\begin{Posledica}
If $T\in\mathcal{B}(\mathcal{H})$ is a finite rank operator, then $\mathcal{R}(T-\lambda)$ is closed for every $\lambda\in\mathds{C}$.
\end{Posledica}

\section{Spectrum}

We begin with results related to the left and the right spectrum, and afterwards conclude with the spectrum.
\begin{Teorema}(\cite[Corollary 2.3]{SARAJLIJA4})\label{POSLEDICA}
Let $D_1\in\mathcal{B}(\mathcal{H}_1),\ D_2\in\mathcal{B}(\mathcal{H}_2),...,D_n\in\mathcal{B}(\mathcal{H}_n)$. Then
\begin{equation}\label{STARIJE4}
\sigma_{l}(D_1)\cup\Big(\bigcup\limits_{k=2}^{n}\Delta_k\Big)\subseteq\\\bigcap\limits_{A\in\mathcal{B}_n}\sigma_{l}(T_n^d(A))
\end{equation}
where
$$
\Delta_k:=\Big\lbrace\lambda\in\mathds{C}:\ \alpha(D_k-\lambda)>\sum\limits_{s=1}^{k-1}\beta(D_s-\lambda)\Big\rbrace,\ 2\leq k\leq n,
$$
\end{Teorema}

\begin{Teorema}\label{SESTA}
Let $D_1\in\mathcal{B}(\mathcal{H}_1),\ D_2\in\mathcal{B}(\mathcal{H}_2),...,D_n\in\mathcal{B}(\mathcal{H}_n)$. Then
\begin{equation}\label{SESTAJEDAN}
\bigcup\limits_{k=1}^n\sigma_{l}(D_k)=\sigma_{l}(T_n^d(A))\cup\Delta
\end{equation}
holds for every $A\in\mathcal{B}_n$, where 
$$
\begin{aligned}
\Delta=\bigcup\limits_{k=2}^{n}\{\lambda\in\sigma_l(D_k):\ \alpha(D_k-\lambda)\leq\sum\limits_{s=1}^{k-1}\beta(D_s-\lambda)\}.
\end{aligned}
$$
\end{Teorema}

\textbf{Proof.} Obviously, $\Delta\subseteq\bigcup\limits_{k=1}^n\sigma_{l}(D_k)$, and $\sigma_l(T_n^d(A))\subseteq\bigcup\limits_{k=1}^n\sigma_{l}(D_k)$ according to Lemma \ref{POMOCNALEMA}. Assume that $\lambda\in\bigcup\limits_{k=1}^n\sigma_{l}(D_k)\setminus\sigma_{l}(T_n^d(A))$. Then by Theorem \ref{POSLEDICA} we get that $\lambda$ does not belong to the left side of (\ref{STARIJE4}), which together with observation $\lambda\in\bigcup\limits_{k=1}^n\sigma_{l}(D_k)$ easily gives $\lambda\in\Delta$. $\square$
\begin{Posledica}
Let $D_1\in\mathcal{B}(\mathcal{H}_1),\ D_2\in\mathcal{B}(\mathcal{H}_2),...,D_n\in\mathcal{B}(\mathcal{H}_n)$. Then
$$
\bigcup\limits_{k=1}^n\sigma_{l}(D_k)=\sigma_{l}(T_n^d(A))
$$
holds for every $A\in\mathcal{B}_n$ if
$$
\Delta=\emptyset,
$$
where $\Delta$ is defined as in Theorem \ref{SESTA}.
\end{Posledica}

If we put $n=2$ we get:

\begin{Teorema}\label{SESTAn=2}
Let $D_1\in\mathcal{B}(\mathcal{H}_1),\ D_2\in\mathcal{B}(\mathcal{H}_2)$. Then
\begin{equation}\label{SESTAJEDAN}
\sigma_{l}(D_1)\cup\sigma_l(D_2)=\sigma_{l}(T_n^d(A))\cup\{\lambda\in\sigma_l(D_2):\ \alpha(D_2-\lambda)\leq\beta(D_1-\lambda)\}
\end{equation}
holds for every $A\in\mathcal{B}_2$.
\end{Teorema}
\begin{Posledica}
Let $D_1\in\mathcal{B}(\mathcal{H}_1),\ D_2\in\mathcal{B}(\mathcal{H}_2)$. Assume that $\mathcal{H}_1,\mathcal{H}_2$ are infinite dimensional Hilbert spaces. Then
$$\sigma_{l}(D_1)\cup\sigma_l(D_2)=\sigma_{l}(T_2^d(A))$$
holds for every $A\in\mathcal{B}_2$ if 
$$\{\lambda\in\sigma_l(D_2):\ \alpha(D_2-\lambda)\leq\beta(D_1-\lambda)\}=\emptyset.$$
\end{Posledica}

Using duality, we obtain results related to the right spectrum.

\begin{Teorema}\label{SEDMA}
Let $D_1\in\mathcal{B}(\mathcal{H}_1),\ D_2\in\mathcal{B}(\mathcal{H}_2),...,D_n\in\mathcal{B}(\mathcal{H}_n)$. Then
\begin{equation}\label{SEDMAJEDAN}
\bigcup\limits_{k=1}^n\sigma_{r}(D_k)=\sigma_{r}(T_n^d(A))\cup\Delta
\end{equation}
holds for every $A\in\mathcal{B}_n$, where 
$$
\begin{aligned}
\Delta=\bigcup\limits_{k=1}^{n-1}\{\lambda\in\sigma_r(D_k):\ \beta(D_k-\lambda)\leq\sum\limits_{s=k+1}^{n}\alpha(D_s-\lambda)\ and\ \mathcal{R}(D_k-\lambda)\ is\ not\ closed \},
\end{aligned}
$$
\end{Teorema}

\begin{Posledica}
Let $D_1\in\mathcal{B}(\mathcal{H}_1),\ D_2\in\mathcal{B}(\mathcal{H}_2),...,D_n\in\mathcal{B}(\mathcal{H}_n)$. Then
$$
\bigcup\limits_{k=1}^n\sigma_{r}(D_k)=\sigma_{r}(T_n^d(A))
$$
holds for every $A\in\mathcal{B}_n$ if
$$
\Delta=\emptyset,
$$
where $\Delta$ is defined as in Theorem \ref{SEDMA}.
\end{Posledica}

Special case $n=2$ gives:

\begin{Teorema}\label{SEDMAn=2}
Let $D_1\in\mathcal{B}(\mathcal{H}_1),\ D_2\in\mathcal{B}(\mathcal{H}_2)$. Then
$$
\sigma_{r}(D_1)\cup\sigma_r(D_2)=\sigma_{r}(T_2^d(A))\cup\{\lambda\in\sigma_r(D_1):\ \beta(D_1-\lambda)\leq\alpha(D_2-\lambda)\}
$$
holds for every $A\in\mathcal{B}_2$.
\end{Teorema}

\begin{Posledica}
Let $D_1\in\mathcal{B}(\mathcal{H}_1),\ D_2\in\mathcal{B}(\mathcal{H}_2)$. Then
$$
\sigma_{r}(D_1)\cup\sigma_r(D_2)=\sigma_{r}(T_2^d(A))
$$
holds for every $A\in\mathcal{B}_2$ if 
$$\{\lambda\in\sigma_r(D_1):\ \beta(D_1-\lambda)\leq\alpha(D_2-\lambda)\}.$$
\end{Posledica}

We finish our investigations with results related to the spectrum of $T_n^d(A)$. First we recall:

\begin{Teorema}(\cite[Corollary 2.14]{SARAJLIJA4})\label{POSLEDICA9}
Let $D_1\in\mathcal{B}(\mathcal{H}_1),\ D_2\in\mathcal{B}(\mathcal{H}_2),...,D_n\in\mathcal{B}(\mathcal{H}_n)$. Then
\begin{equation}\label{STARIJE5}
\sigma_{l}(D_1)\cup\sigma_{r}(D_n)\cup\Big(\bigcup\limits_{k=2}^{n-1}\delta_k\Big)\cup\delta_n\subseteq\\\bigcap_{A\in\mathcal{B}_n}\sigma(T_n^d(A)),
\end{equation}
where
$$
\begin{aligned}
\delta_k=\Big\lbrace\lambda\in\mathds{C}:\ \alpha(D_k-\lambda)>\sum_{s=1}^{k-1}\beta(D_s-\lambda)\Big\rbrace\cup\\
\Big\lbrace\lambda\in\mathds{C}:\ \beta(D_k-\lambda)>\sum_{s=k+1}^n\alpha(D_s-\lambda)\Big\rbrace,\quad 2\leq k\leq n-1,
\end{aligned}
$$
$$
\begin{aligned}
\delta_n=\Big\lbrace\lambda\in\mathds{C}:\ \alpha(D_n-\lambda)>\sum_{s=1}^{n-1}\beta(D_s-\lambda)\Big\rbrace\cup\\\Big\lbrace\lambda\in\mathds{C}:\ \beta(D_1-\lambda)>\sum_{s=2}^n\alpha(D_s-\lambda)\Big\rbrace.
\end{aligned}
$$
\end{Teorema}

\begin{Teorema}\label{OSMA}
Let $D_1\in\mathcal{B}(\mathcal{H}_1),\ D_2\in\mathcal{B}(\mathcal{H}_2),...,D_n\in\mathcal{B}(\mathcal{H}_n)$. Then
$$
\bigcup\limits_{k=1}^n\sigma(D_k)=\sigma(T_n^d(A))\cup\Delta
$$
holds for every $A\in\mathcal{B}_n$, where 
$$
\begin{aligned}
\Delta=\rho_l(D_1)\cap\rho_r(D_n)\cap\bigcup\limits_{k=2}^{n-1}\Big\{\lambda\in\mathds{C}:\ \Big(\beta(D_k-\lambda)\leq\sum\limits_{s=k+1}^{n}\alpha(D_s-\lambda)\ and\\ \lambda\in\sigma_r(D_k)\Big)\ or\ \Big(\alpha(D_k-\lambda)\leq\ \sum\limits_{s=1}^{k-1}\beta(D_s-\lambda\Big)\ and\ \lambda\in\sigma_l(D_k)\Big)\Big\}\bigcup\\ \Big\{\lambda\in\mathds{C}:\ 0<\beta(D_1-\lambda)\leq\sum\limits_{s=2}^{n}\alpha(D_s-\lambda)\ or\ 0<\alpha(D_n-\lambda)\leq\sum_{s=1}^{n-1}\beta(D_s-\lambda)\Big\}.
\end{aligned}
$$
\end{Teorema}
\textbf{Proof.} First of all, $\sigma(T_n^d(A))\subseteq\bigcup\limits_{k=1}^n\sigma(D_k)$ according to Lemma \ref{POMOCNALEMA}. Next, we prove that $\Delta\subseteq\bigcup\limits_{k=1}^n\sigma(D_k)$. This is, however, obvious from the definition of an invertible operator. Namely, if $\lambda\in\Delta$, then either we have $\lambda\in\sigma_r(D_k)\cup\sigma_l(D_k)$ for some $k\in\lbrace2,...,n-1\rbrace$ or $\beta^2(D_1-\lambda)+\alpha^2(D_n-\lambda)>0$. In both cases we get $\lambda\in\bigcup\limits_{k=1}^n\sigma(D_k)$. Now, let us prove the opposite inclusion. Assume that $\lambda\in\bigcup\limits_{k=1}^n\sigma(D_k)\setminus\sigma(T_n^d(A))$. However, by Theorem \ref{POSLEDICA9}, we get that $\lambda$ does not belong to the left side of (\ref{STARIJE5}), which using the notation from (\ref{STARIJE5}), yields $\lambda\in\rho_l(D_1)\cap\rho_r(D_n)\cap\bigcap\limits_{k=2}^n\delta_k^c$. This, together with our assumption $\lambda\in\bigcup\limits_{k=1}^n\sigma(D_k)$ easily gives $\lambda\in\Delta$. $\square$
\begin{Posledica}
Let $D_1\in\mathcal{B}(\mathcal{H}_1),\ D_2\in\mathcal{B}(\mathcal{H}_2),...,D_n\in\mathcal{B}(\mathcal{H}_n)$. Then
$$
\bigcup\limits_{k=1}^n\sigma(D_k)=\sigma(T_n^d(A))
$$
holds for every $A\in\mathcal{B}_n$ if
$$
\Delta=\emptyset,
$$
where $\Delta$ is defined as in Theorem \ref{OSMA}.
\end{Posledica}

Statements related to the spectrum become especially elegant when $n=2$.

\begin{Teorema}\label{OPSHTIJE}
Let $D_1\in\mathcal{B}(\mathcal{H}_1),\ D_2\in\mathcal{B}(\mathcal{H}_2)$. Then
$$
\begin{aligned}
\sigma(D_1)\cup\sigma(D_2)=\sigma(T_2^d(A))\cup\Big\{\lambda\in\mathds{C}:\ 0<\beta(D_1-\lambda)\leq\alpha(D_2-\lambda)\ or\\ 0<\alpha(D_2-\lambda)\leq\beta(D_1-\lambda)\Big\}.
\end{aligned}
$$
holds for every $A\in\mathcal{B}_2$.
\end{Teorema}
Observe that Theorem \ref{OPSHTIJE} is a more precise version of \cite[Corollary 7]{HAN} in the Hilbert space setting. Namely, in \cite{HAN}, authors prove that a passage from $\sigma(D_1)\cup\sigma(D_2)$ to $\sigma(T_2^d(A))$ is accomplished by filling some holes in $\sigma(T_2^d(A))$ which happen to be subsets of $\sigma(D_1)\cap\sigma(D_2)$. Notice, however, that in Theorem \ref{OPSHTIJE}, we have specified the form of these holes. To our best knowledge, this has not been done so far. 
\begin{Posledica}
Let $D_1\in\mathcal{B}(\mathcal{H}_1),\ D_2\in\mathcal{B}(\mathcal{H}_2)$. Then
$$
\sigma(D_1)\cup\sigma(D_2)=\sigma(T_2^d(A))
$$
holds for every $A\in\mathcal{B}_2$ if
$$
\begin{aligned}
\Big\{\lambda\in\mathds{C}:\ 0<\beta(D_1-\lambda)\leq\alpha(D_2-\lambda)\ or\\ 0<\alpha(D_2-\lambda)\leq\beta(D_1-\lambda)\Big\}=\emptyset.
\end{aligned}
$$
\end{Posledica}

\noindent\textbf{Statements \& Declarations}\\[3mm]
\hspace*{6mm}This work was supported by the Ministry of Education, Science and Technological Development of the Republic of Serbia under Grant No. 451-03-66/2024-03/200125 and 451-03-65/2024-03/200125.\\[1mm]
\hspace*{6mm}The author has no relevant financial or non-financial interests to disclose.

\noindent\textbf{Data availability statement}\\[3mm]
\hspace*{6mm}Data sharing not applicable to this article as no datasets were generated or analysed during the current study.
\end{document}